\documentclass[9pt]{article}
\usepackage{amsfonts}
\usepackage{mathrsfs}
\usepackage{ifpdf}
\usepackage{amsfonts, amsmath, amssymb}
\usepackage{amssymb,amsfonts,amsmath,
latexsym, epsfig,cite, psfrag,eepic,colordvi}
\usepackage{amscd,graphics}
\newtheorem{theorem}{Theorem}[section]

\newtheorem{lemma}[theorem]{Lemma}
\newtheorem{corollary}[theorem]{Corollary}
\numberwithin{equation}{section}

\def\qed{\nopagebreak\hfill{\rule{4pt}{7pt}}\medbreak}
\begin{document}

\begin{center}
{\large {ON FLUSHED PARTITIONS AND CONCAVE COMPOSITIONS}}

\vskip 6mm {\small Xiao-Chuan Liu\footnote{The project is supported
partially by CNNSF
(No.10971106).}\\[%
2mm] Nankai University, Tianjin 300071,
P.R. China \\[3mm]
E-mail: lxc1984@gmail.com \\[0pt%
] }
\end{center}

\noindent{\bf Abstract.}\ In this work, we give combinatorial proofs
for generating functions of two problems, i.e., flushed partitions
and concave compositions of even length. We also give combinatorial
interpretation of one problem posed by Sylvester involving flushed
partitions and then prove it. For these purposes, we first describe
an involution and use it to prove core identities. Using this
involution with modifications, we prove several problems of
different nature, including Andrews' partition identities involving
initial repetitions and partition theoretical interpretations of three
mock theta functions of third order $f(q)$, $\phi(q)$ and $\psi(q)$.
An identity of Ramanujan is proved combinatorially. Several new
identities are also established.



\section{Introduction}

In this paper, we are mainly concerned with two problems in the
theory of integer partition, namely, flushed partitions and concave
compositions of even length.

The definition of flushed partition is given by Sylvester
\cite{Sylvester}. A partition is called {\it flushed} when the
number of the parts with length $k$ is odd, where
$k=1,2,\cdots,2i-1$, and the parts with length $2i$ do not
occur an odd number of times. Similarly, define {\it unflushed
partitions} as those not satisfying the above conditions. Sylvester
also posed two problems with respect to flushed partitions. One of
them, in Sylvester's words, is stated as follows,

`` 1. Required to prove, that if any number be partitioned in every
possible way, the number of unflushed partitions containing an odd
number of parts is equal to the number of unflushed partitions
containing an even number of parts.

`` Ex.gr.: The total partitions of $7$ are $7;$ $6, 1;$ $5,2;$ $5,
1, 1;$ $4, 3;$ $4, 2, 1;$ $4, 1, 1, 1;$ $3, 3, 1;$ $3, 2, 2;$ $3, 2,
1, 1;$ $2, 2, 2, 1;$ $3, 1, 1, 1, 1;$ $2, 2, 1, 1, 1;$ $2, 1, 1, 1,
1, 1;$ $1, 1, 1, 1, 1, 1, 1.$ Of these, $6, 1;$ $4, 1, 1, 1;$ $3, 3,
1;$ $2, 2, 1, 1, 1;$ $1, 1, 1, 1, 1, 1, 1$ alone are flushed. Of the
remaining unflushed partitions, five contain an odd number of parts,
and five an even number.

``Again, the total partitions of $6$ are $6;$ $5, 1;$ $4, 2;$ $4, 1,
1;$ $3, 3;$ $3, 2, 1;$ $2, 2, 2;$ $3, 1, 1, 1;$ $2, 2, 1, 1;$ $2, 1,
1, 1, 1;$ $1, 1, 1, 1, 1, 1;$ of which $5, 1;$ $3, 2, 1;$ $3, 1, 1,
1$ alone are flushed. Of the remainder, four contain an odd and four
an even number of parts.

``N.B.|This transcendental theorem compares singularly with the
well-known algebraical one, that the total number of the permuted
partitions of a number with an odd number of parts is equal to the
same of the same with an even number.

Solution to this problem is given by Andrews in 1970 in
\cite{GeorgeAndrews5} by manipulating generating functions, and the
generating function for flushed partitions reads as follows,

\begin{equation}\frac{1}{(q)_\infty}\sum_{n=1}^\infty q^{n(3n-1)/2}(1-q^n).\end{equation}

 However, Andrews doubted that his proofs are what
Sylvester expected in the first place. He writes (in \cite{Andrews7}):
``It is completely unknown whether this was
Sylvester's approach and how he came upon flushed partitions in the
first place."

It is hardly believed that the core of the combinatorial proof of
Sylvester's problem is only one involution. But it turns out to be
the case. We will first prove combinatorially the above generating
function (1.1), and then by inserting another variable $z$ into the
generating function, we can map bijectively unflushed partition of
$n$ with $m$ parts into unrestricted partitions of $n$ with $m$ parts
with additional restrictions involving Durfee Symbols, a very
natural concept yet just introduced recently by Andrews in
\cite{GeorgeAndrews1}. In this way, we try to understand flushed
partitions in a new combinatorial sense and then provide a new proof
of Sylvester's problem.

{\it Concave composition of even length} was recently introduced by
Andrews in the study of orthogonal polynomials, see \cite{GeorgeAndrews2, GeorgeAndrews3}.
It is a sum of the form $\sum
a_i+\sum b_i$ such that
$$a_1>a_2>\cdots>a_m=b_m<b_{m-1}<\cdots<b_1,$$
where $a_m\geq 0$, and all $a_i$ and $b_i$ are integers. Let
$\mathcal {C}\mathcal {E}(n)$ denote the set of concave compositions
of even length of $n$, and let $ce(n)$ be the cardinality of
$\mathcal {C}\mathcal {E}(n)$. By transformation formulas, Andrews
derived the generating function of $ce(n)$ as follows
([\mdseries{7, Theorem 1}]),

For $|q|<1$,
\begin{equation}\label{1}
 \sum_{n=0}^\infty
ce(n)q^n=\frac{1}{(q)_\infty} \Big(1-\sum_{n=1}^\infty
    q^{n(3n-1)/2}(1-q^n)\Big).
\end{equation}

Andrews \cite{GeorgeAndrews2} asked for a combinatorial proof of
Theorem $\ref{1}$. We will give one such proof in this paper.

Note the above two generating functions (1.1) and (1.2) have
connections with one another, which lead to our main theorem of this
paper as stated as follows.

\begin{theorem}\label{11} The number of unflushed partitions of n is equal to the number
of concave compositions of even length of $n$.
\end{theorem}

This paper is organized as follows. In Section $2$, we define an involution,
which, with modifications, will be used
repeatedly. In Section $3$, we give combinatorial proofs of several
different q-series identities. Readers who are not interested in
these problems can skip directly to Section $4$, where we will prove
our main theorem about these two generating functions and prove
Sylvester's first problem.

The applications of the involutions in Section $3$ consists of
several different partition identities. The nature of these problems
varies, showing such involutions are indeed useful tools.

One application is on Andrews' partition identity involving initial
repetitions. In \cite{Andrews6}, Andrews proved the q-series
identity

\begin{equation}
    \sum_{n=0}^{\infty}\frac{z^nq^{1\cdot 2+2\cdot 2+\cdots +n\cdot
2}}{(q)_n}\prod_{j=n+1}^{\infty}(1-zq^j)= \sum_{j=0}^\infty
(-1)^jz^jq^{\frac{j(j+1)}{2}},
\end{equation}
which, when interpreted combinatorially, means that partitions of
$n$ with $m$ different parts and an even number of distinct parts in
which, if part $j$ is repeated, then all parts smaller than $j$ is
repeated, are equinumerous with partitions of $n$ with $m$ different
parts and an odd number of distinct parts satisfying the same
property, unless $n$ is a triangular number $\frac{j(j+1)}{2}$ and
$m=j$, when their difference is $(-1)^j$. Partitions with this
property are called {\it partitions with initial 2-repetitions}.

This result (1.3) looks similar to classical Euler's pentagonal
number theorem. The involution plays a role as the role the Franklin's
well known involution has played in Euler's pentagonal number
theorem. Based on this proof, a formula is given to compute the
number of partitions into even number of distinct parts.

The next applications involve combinatorial interpretations of the
following three mock theta functions of order $3$, defined by
Ramanujan in his last letter to Hardy (see \cite{lastletter}),

\begin{equation}\label{f14}
f(q)=\sum_{n=0}^\infty\frac{q^{n^2}}{(-q)_n^2}
\end{equation}
\begin{equation}\label{f15}
\phi(q)=\sum_{n=0}^\infty\frac{q^{n^2}}{(-q^2;q^2)_n}
\end{equation}
\begin{equation}\label{f16}
\psi(q)=\sum_{n=1}^\infty\frac{q^{n^2}}{(q;q^2)_n}
\end{equation}

The combinatorial interpretation of (\ref{f14}) and (\ref{f15}) are
given in the end of Section $2$. The following combinatorial
interpretation of (\ref{f16}) was first given by Fine in \cite{Fine2}:
\begin{equation}f(q)=1+\frac{1}{(-q)_\infty}\sum_{k\geq 1}
(-1)^{k-1}q^k(-q^{k+1})_\infty
\end{equation}
We will restate this result as Theorem \ref{37}. Our proof is new.

Then all these combinatorial interpretations lead to the following
identity of Ramanujan:

$$f(q)=\phi(-q)-2\psi(-q)$$

We will restate this as Corollary 3.8. The first proof of this
identity is given by Watson in \cite{Watson}. Another combinatorial
proof was given by Chen, Ji and Liu in \cite{chen}.

At the end of Section $3$, we use the involution to generate several
more identities, which all have partition interpretations. We only
write such interpretation for one identity as an example.

\section{An Involution $\alpha$}
All notations in the theory of integer partitions follow the book
\cite{GeorgeAndrews1}. A partition of a positive integer $n$
is a finite nonincreasing sequence of positive integers
$(\lambda_1,\cdots,\lambda_r)$ such that $\sum_{i=1}^r \lambda_i=n$.
$\lambda_i$ are called the parts of the partition. We use $\mathcal
P$ to denote the set of partitions and $\mathcal P_n$ to denote the
set of partitions of $n$. We are also interested in the partitions into
distinct parts. We use $\mathcal D$ to denote the set of partitions into
distinct parts and $\mathcal D_n$ to denote the set of partitions of
$n$ into distinct parts.

At the same time, define the set of partitions into distinct parts which may contain
one copy of empty part as $\mathcal D'$. Naturally, $\mathcal D'_n$
would denote the set of such partitions of $n$. Similarly, we define
$\mathcal P'$ as the set of partitions which can contain empty parts
(maybe more than one yet a limited number of copies). $\mathcal P'_n$, of course, would
be the set of such partitions of $n$.

We can represent a partition as its
Young diagram, that is, a pattern of left-justified boxes with
$\lambda_i$ squares in row $i$. The square in the $i$th row and the
$j$th column can be written simply as the square $(i,j)$.  The
Durfee square in $\lambda$ is the largest square of boxes contained
in the partition $\lambda$. When the Durfee square of a partition is
$n\times n$, we say it is a Durfee square of size $n$. The conjugate
of $\lambda=(\lambda_1,\cdots,\lambda_r)$ is a partition $\lambda'$
with the $i$th part $\lambda_i'$ as the number of parts of $\lambda$
that are $\geq i$.

We adopt the following standard notation:
 \begin{align*}
 &(a;q^k)_n=(1-a)(1-aq^k)\cdots(1-aq^{(n-1)k}),\\
 &(a;q^k)_\infty=\prod_{n=0}^\infty(1-aq^{nk}),\\
 &(a)_n=(a;q)_n,\\
 &(a)_0=(a;q)_0=1.\end{align*}

Now we are ready to state the involution we mentioned, which will be denoted as
$\alpha$. Given a triple $(\lambda,\mu,\rho_{n+\ell})$ satisfying
the following properties with a sign $(-1)^{k+\ell}$.

1. $\lambda\in\mathcal D'$ has n distinct parts. The least part has length $k$. Note that
it may contain one empty part.

2. $\mu$ is a partition with no more than $\ell$ parts. Or, we can think
$\mu\in \mathcal P$ as a partition into exactly $\ell$ parts, which may contain
some empty parts.

3. $\rho_{n+\ell}$ is Sylvester's triangle
$(n+\ell,n+\ell-1,\cdots,1)$.

Let the involution $\alpha$ act on such a triple by comparing
$\lambda_1$ the largest part of $\lambda$ and $\mu_1$ the largest
part of $\mu$, as follows.

If $\lambda_1\geq \mu_1$, and the number of parts of $\lambda$ is at least two,
 then we move the first part of $\lambda$
and attach it to $\mu$, making it $\mu'$. What has been left there
is then $\lambda'$. Since $\mu'$ has no more than $\ell+1$ parts,
yet the smallest part of $\lambda$ does not change, the sign changes.

On the other hand, if $\lambda_1< \mu_1$, then we move the first
part of $\mu$ and attach it to $\lambda$, making it $\lambda'$. What
has been left there is then $\mu'$ and $\mu'$ has no more than
$\ell-1$ parts, yet the smallest part of $\lambda$ does not change,
the sign changes too.

The triples that the involution does not apply to are the ones
$(\lambda,\mu,\rho_d)$ where $\lambda$ has only one part of length
$t$ and $t\geq \mu_1$; $\mu$ is a partition with no more than $d-1$
parts.

\begin{figure}[h,t] \setlength{\unitlength}{0.4mm}
\begin{center}
  \begin{picture}(350,170)
    \put(0,0) {\begin{picture}(110,110)
      \put(0,120){\line(1,0){20}}
      \put(0,120){\line(0,1){40}}

      \put(0,160){\line(1,0){90}}
      \put(0,150){\line(1,0){90}}
 \put(90,160){\line(0,-1){10}}
  \put(90,150){\line(-1,0){30}}
  \put(60,140){\line(-1,0){10}}
  \put(60,140){\line(0,1){10}}
  \put(50,130){\line(0,1){10}}
  \put(50,130){\line(-1,0){30}}
    \put(20,120){\line(0,1){10}}

\put(60,10){\line(0,1){70}}

\thicklines

       \put(60,10){\line(0,1){40}}
       \put(0,80){\line(1,0){130}}
       \put(60,10){\line(1,0){10}}
       \put(70,20){\line(1,0){10}}
       \put(70,20){\line(0,-1){10}}
       \put(80,30){\line(1,0){10}}
       \put(80,30){\line(0,-1){10}}
       \put(90,40){\line(1,0){10}}
       \put(90,40){\line(0,-1){10}}
       \put(100,50){\line(1,0){10}}
       \put(100,50){\line(0,-1){10}}
       \put(110,60){\line(1,0){10}}
       \put(110,60){\line(0,-1){10}}
       \put(120,60){\line(0,1){10}}
       \put(130,70){\line(-1,0){10}}
       \put(130,70){\line(0,1){10}}

       \put(0,70){\line(0,1){10}}
       \put(0,70){\line(1,0){20}}
       \put(20,50){\line(0,1){20}}
       \put(20,50){\line(1,0){40}}
       \thinlines
\put(0,100){$ \text{ the partition }\lambda=(9,6,5,2)$} \put(0,0){$
\text{ the partition }\mu=(6,4,4) \text { with }  \rho_7$}
\end{picture}
} \thicklines
   \put(150,90){\vector(1,0){25}}
  \thinlines
      \put(220,0){\begin{picture}(110,100)
      \put(0,120){\line(1,0){20}}
      \put(0,120){\line(0,1){30}}
      \put(0,150){\line(1,0){60}}

      \put(60,140){\line(-1,0){10}}
      \put(60,140){\line(0,1){10}}
      \put(50,130){\line(0,1){10}}
      \put(50,130){\line(-1,0){30}}
      \put(20,120){\line(0,1){10}}

      \put(60,10){\line(0,1){70}}
      \put(-30,80){\line(1,0){160}}
      \put(-30,80){\line(0,-1){10}}
      \put(-30,70){\line(1,0){90}}
\thicklines

       \put(60,10){\line(0,1){30}}

       \put(60,10){\line(1,0){10}}
       \put(70,20){\line(1,0){10}}
       \put(70,20){\line(0,-1){10}}
       \put(80,30){\line(1,0){10}}
       \put(80,30){\line(0,-1){10}}
       \put(90,40){\line(1,0){10}}
       \put(90,40){\line(0,-1){10}}
       \put(100,50){\line(1,0){10}}
       \put(100,50){\line(0,-1){10}}
       \put(110,60){\line(1,0){10}}
       \put(110,60){\line(0,-1){10}}
       \put(120,60){\line(0,1){10}}
       \put(130,70){\line(-1,0){10}}
       \put(130,70){\line(0,1){10}}

       \put(0,70){\line(1,0){60}}
       \put(0,60){\line(0,1){10}}
       \put(0,60){\line(1,0){20}}
       \put(20,40){\line(0,1){20}}
       \put(20,40){\line(1,0){40}}
       \thinlines
\put(0,100){$ \text{ the partition }\lambda'=(6,5,2)$} \put(0,0){$
\text{ the partition }\mu'=(9,6,4,4) \text { with }  \rho_7$}
\end{picture}}
\end{picture}\caption{ The involution $\alpha$}

\end{center}
\end{figure}
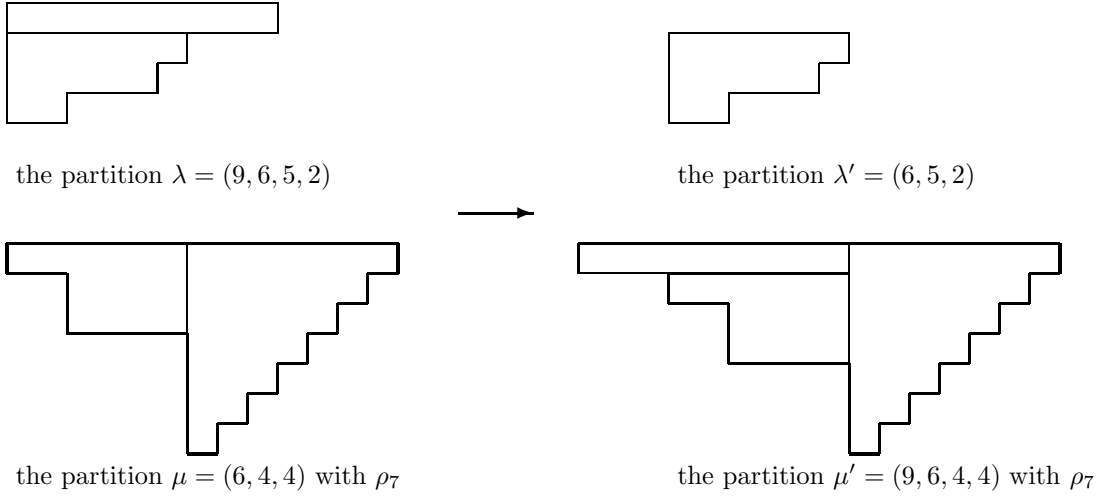

For example, Let $\lambda=(9,6,5,2),\mu=(6,4,4),k=2,n=4,\ell=3$. It
is assigned $(-1)^5=-1$. Then we have that $\alpha
(\lambda,\mu,\rho_{4+3})=(\lambda',\mu',\rho'_{3+4})$, where
$\lambda'=(6,5,2),\mu'=(9,6,4,4),k=2,n'=3,\ell'=4$. It is assigned
with $(-1)^6=1$. We illustrate this example in Figure $1$.

\noindent\textbf{Remark 1.} For the triple
$(\lambda,\mu,\rho_{n+\ell})$, special attention should be given to $\mu$, because
sometimes the number of nonzero parts of $\mu$ is strictly less than $\ell$. This situation makes no
exception for our arguments. In fact, $(\lambda',\mu',\rho_{(n\pm
1)+(\ell\mp 1)})=\alpha ((\lambda,\mu,\rho_{n+\ell}))$ where $\ell$
changes by $\pm 1$.

\noindent\textbf{Remark 2.} The nature of this involution has a more
general form in \cite{Pak}.

Applying this involution, we can prove the following identity:

\begin{theorem}
\begin{equation}\label{21}
\sum_{n=1}^\infty\frac{q^{n^2}(q^{n+1})_\infty}{(q)_{n-1}(1+q^n)}
 =\sum_{n=1}^\infty\frac{(-1)^{n-1}q^{\frac{n(n+1)}{2}}}{(1+q)\cdots(1+q^n)}.
\end{equation}
\end{theorem}
\vskip 0.6cm

\noindent\textbf{Proof.}

For each term in the left hand side of (\ref{21}), we interpret
$q^{n^2},\frac{1}{(q)_{n-1}},\frac{1}{1+q^n} \text{ and }
(q^{n+1})_\infty$ as follows, respectively:

\noindent 1. A $n\times n$ squares,

\noindent 2. A partition with the largest part at most
$n-1$,

\noindent 3. $k$ parts of length exactly $n$, assigned $(-1)^k$,

\noindent 4. A partition with distinct parts, the least part being
larger than $n$. (Suppose this partition has $\ell$ parts. It is
assigned $(-1)^\ell$.)

We split the $n\times n$ squares to a Sylvester's triangle
$\rho_n=(n,n-1,\cdots,1)$ and another Sylvester's triangle
$\rho_{n-1}=(n-1,n-2,\cdots,1)$. Then we glue together three objects
to form a new partition $\lambda^*$, which are, the $k$ parts of
length $n$, $\rho_{n-1}$ and the partition with the largest part at most $n-1$.

Now we observe its conjugate and denote it as $\lambda$. It is a
partition with distinct parts, $k$ the length of the least part.
Obviously, $k=0$ is allowed, so $\lambda$ can have
one empty part, that is, $\lambda\in \mathcal D'$.

Then we attach the Sylvester's triangle $\rho_{n}$ under the
partition with $\ell$ parts, the least length of parts larger than
$n$. This is clearly a partition with distinct parts. We divide it
into a partition $\mu$ and Sylvester's triangle $\rho_{n+\ell}$ in
the obvious way. We see that $\mu$ has no more than $\ell$ parts.
Again, we can regard $\mu$ as a partition with exactly $\ell$ parts, where
we allow empty parts to exist.

Remember this was assigned $(-1)^{k+\ell}$.

We then invoke the involution $\alpha$. Since the cases which the involution
does not apply to are the triples $(\lambda,\mu,\rho_d)$
where $\lambda$ has only one part of length $t$ and $t\geq \mu_1$;
$\mu$ is a partition with no more than $d-1$ parts.

We move this $t$ squares and attach it into $\mu$ together, getting
a partition with no more than $d$ parts, which is assigned $(-1)^{t+d-1}$,
$t$ the largest part of $\mu$, $d$ is the subscript
of the Sylvester's Triangle. This corresponds to the right-hand side of
the identity, which completes the proof. \qed

\begin{figure}[h,t]
\setlength{\unitlength}{0.4mm}
\begin{center}
  \begin{picture}(300,170)
    \put(0,0) {\begin{picture}(135,135)

    \put(0,120){\line(1,0){130}}

    \put(0,110){\line(1,0){100}}
    \put(0,100){\line(1,0){90}}

    \put(0,80){\line(1,0){40}}
    \put(0,40){\line(1,0){40}}
    \put(0,0){\line(1,0){10}}

    \put(0,0){\line(0,1){130}}
    \put(10,0){\line(0,1){20}}
    \put(10,20){\line(1,0){20}}
    \put(30,20){\line(0,1){20}}
 \put(40,80){\line(0,1){20}}
\thicklines
    \put(0,100){\line(1,0){40}}
    \put(0,100){\line(0,1){30}}
    \put(0,80){\line(1,0){40}}
    \put(0,70){\line(1,0){10}}
    \put(0,70){\line(0,1){10}}
    \put(10,60){\line(1,0){10}}\put(10,60){\line(0,1){10}}
    \put(20,50){\line(1,0){10}}\put(20,50){\line(0,1){10}}
    \put(30,40){\line(0,1){10}} \put(30,40){\line(1,0){10}}
    \put(40,40){\line(0,1){40}}\put(40,100){\line(1,0){50}}
     \put(130,120){\line(0,1){10}}
    \put(100,110){\line(0,1){10}}
    \put(90,100){\line(0,1){10}}
    \put(90,110){\line(1,0){10}}
    \put(100,120){\line(1,0){30}}
    \put(0,130){\line(1,0){130}}
    \thinlines
\end{picture}}
\thicklines
   \put(130,50){\vector(1,0){30}}
  \thinlines
      \put(180,0){\begin{picture}(110,100)
      \put(0,120){\line(1,0){20}}
      \put(0,120){\line(0,1){40}}
      \put(0,160){\line(1,0){90}}
 \put(90,160){\line(0,-1){10}}
  \put(90,150){\line(-1,0){30}}
  \put(60,140){\line(-1,0){10}}
  \put(60,140){\line(0,1){10}}
  \put(50,130){\line(0,1){10}}
  \put(50,130){\line(-1,0){30}}
    \put(20,120){\line(0,1){10}}

\put(60,10){\line(0,1){70}}

\thicklines

       \put(60,10){\line(0,1){40}}
       \put(0,80){\line(1,0){130}}
       \put(60,10){\line(1,0){10}}
       \put(70,20){\line(1,0){10}}
       \put(70,20){\line(0,-1){10}}
       \put(80,30){\line(1,0){10}}
       \put(80,30){\line(0,-1){10}}
       \put(90,40){\line(1,0){10}}
       \put(90,40){\line(0,-1){10}}
       \put(100,50){\line(1,0){10}}
       \put(100,50){\line(0,-1){10}}
       \put(110,60){\line(1,0){10}}
       \put(110,60){\line(0,-1){10}}
       \put(120,60){\line(0,1){10}}
       \put(130,70){\line(-1,0){10}}
       \put(130,70){\line(0,1){10}}

       \put(0,70){\line(0,1){10}}
       \put(0,70){\line(1,0){20}}
       \put(20,50){\line(0,1){20}}
       \put(20,50){\line(1,0){40}}
       \thinlines
\put(0,100){$ \text{ the partition }\lambda=(9,6,5,2)$} \put(0,0){$
\text{ the partition }\mu=(6,4,4) \text { with } \rho_7$}
\end{picture}}
\end{picture}\caption{ representation of left-hand side of $(2.1)$}

\end{center}
\end{figure}
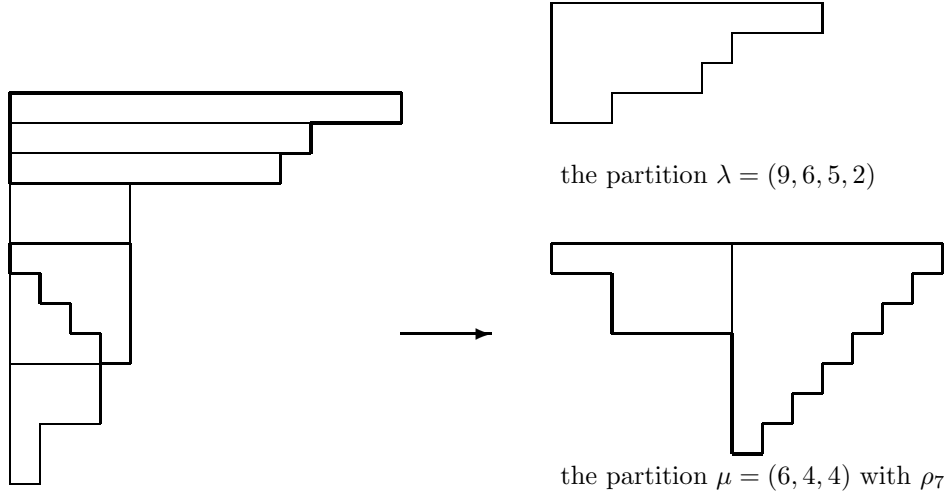

See Figure 2 for a concrete example. The graph in the left is the
partition $(13, 10, 9, 4, 4, 4, 4, 4,$ $4, 3, 3, 1, 1)$, while the
graph in the right is $\lambda=(9,6,5,2), \mu=(6,4,4) \text{ and }
\rho_{4+3}$.

Now we generalize theorem  \ref{21} to the following two identities:

\begin{corollary}

\begin{equation}\label{f22}
\sum_{n=1}^\infty\frac{z^nq^{n^2}(zq^{n+1})_\infty}{(q)_{n-1}(1+q^n)}
 =\sum_{n=1}^\infty\frac{(-1)^{n-1}z^nq^{\frac{n(n+1)}{2}}}{(1+q)\cdots(1+q^n)}.
\end{equation}

\begin{equation}\label{f23}
\sum_{n=1}^\infty\frac{z^nq^{n^2}(zq^{n+1})_\infty}{(q)_{n-1}(1+zq^n)}
 =\sum_{n=1}^\infty\frac{(-1)^{n-1}z^nq^{\frac{n(n+1)}{2}}}{(1+zq)\cdots(1+zq^n)}.
\end{equation}
\end{corollary}

\noindent\textbf{Proof.}
For the first identity, recall the first remark we made after the definition of
the involution $\alpha$. We observe that after the application of the involution $\alpha$,
the size of the Sylvester's triangle does not change. We use another variable $z$ to track
this size. Following the proof of Theorem \ref{21}, we write $z^{n+\ell}$ to denote the Sylvester's
Triangle $\rho_{n+\ell}$. In this way, we have proved the identity (\ref{f22}).

For the second identity, we insert the variable $z$ to track the sum of the size of the Sylvester's triangle
and the number of additional parts of length $n$ apart from the $n\times n$ squares. In the context
of the proof of theorem \ref{21}, we write $z^{n+\ell+k}$ to denote this value. Again, since this value
does not change under the involution $\alpha$, we have proved the identity (\ref{f23}).
\qed

In (\ref{f22}) and in (\ref{f23}), let $q\to q^2, z\to -q^{-1}$, we
have the following:

\begin{equation}\label{f24}
\sum_{n=1}^\infty\frac{(-1)^nq^{2n^2-n}(-q^{2n+1};q^2)_\infty}{(q^2;q^2)_{n-1}(1+q^{2n})}
 =-\sum_{n=1}^\infty\frac{q^{n^2} } {(-q^2;q^2)_n}=1-\phi(q),
\end{equation}

\begin{equation}\label{f25}
\sum_{n=1}^\infty\frac{(-1)^nq^{2n^2-n}(-q^{2n+1};q^2)_\infty}{(q^2;q^2)_{n-1}(1-q^{2n-1})}
 =-\sum_{n=1}^\infty\frac{q^{n^2} } {(q;q^2)_n}=-\psi(q),
\end{equation}
where
$$ \phi(q)=\sum_{n=0}^\infty\frac{q^{n^2}}{(-q^2;q^2)_n}, $$
$$\psi(q)=\sum_{n=1}^\infty\frac{q^{n^2}}{(q;q^2)_n},$$
are two mock theta functions of third order defined by Ramanujan (see \cite{lastletter}).

By simple calculations with Euler's identity
$(-q)_\infty=\frac{1}{(q;q^2)_\infty}$, we have

\begin{equation}\label{f26}
(-q;q)_\infty\big (1-\phi(-q)\big)=
\sum_{n=1}^\infty\frac{q^{2n^2-n}}{(q)_{2n-1}(1+q^{2n})},
\end{equation}

\begin{equation}\label{f27}
(-q;q)_\infty \big(-\psi(-q)\big)=
\sum_{n=1}^\infty\frac{q^{2n^2-n}}{(q)_{2n-1}(1+q^{2n-1})}=
\sum_{n=1}^\infty\frac{q^{2n^2-n}}{(q)_{2n-2}(1-q^{2(2n-1)})}.
\end{equation}

We now use $\mathcal D_{e,e}(n)$ to denote the number of partitions
of $n$ into an even number of distinct parts in which the smallest part
is even. $\mathcal D_{e,o}(n)$ denote the number of partitions of
$n$ into an even number of distinct parts in which the smallest part is
odd.

Symmetrically, $\mathcal D_{o,e}(n)$ (respectively, $\mathcal
D_{o,o}(n)$) denote the number of partitions of
$n$ into odd number of distinct parts in which the smallest part is
even (respectively, odd).

Then by interpreting (\ref{f26}) and (\ref{f27}) we get the
partition theoretical interpretations of mock theta function $\phi(-q)$
and $\psi(-q)$:

\begin{theorem}
\begin{equation}\label{f28}
(-q;q)_\infty \big(1-\phi(-q)\big)= \sum_{n=1}^\infty(\mathcal
D_{o,e}(n)+\mathcal D_{o,o}(n)+\mathcal D_{e,e}(n)-\mathcal
D_{e,o}(n))q^n,\end{equation}
\begin{equation}\label{f29}
(-q;q)_\infty \big(-\psi(-q)\big)=\sum_{m=1}^\infty \mathcal
D_{o,o}(n)q^n.
\end{equation}

\end{theorem}

\noindent\textbf{Proof.}
For identity (\ref{f28}), note that $q^{2n^2-n}=q^{1+2+\cdots+2n-1}$, and that $\frac{1}{(1+q^{2n})}=\sum_{k=0}^\infty
(-q^{2n})^k$. Then the right-hand side of the identity (\ref{f26}) is the sum of the following three parts,

$$\sum_{n=1}^\infty\frac{q^{2n^2-n}}{(q)_{2n-1}},\ \sum_{n=1}^\infty\frac{q^{2n^2-n}}{(q)_{2n-1}}(q^{4n}+q^{8n}+\cdots) \text{ and }-\sum_{n=1}^\infty\frac{q^{2n^2-n}}{(q)_{2n-1}}(q^{2n}+q^{6n}+\cdots),$$
which corresponds to
$$\sum_{n=1}^\infty(\mathcal
D_{o,e}(n)+\mathcal D_{o,o}(n)q^n,\ \sum_{n=1}^\infty \mathcal D_{e,e}(n)q^n \text{ and } -\sum_{n=1}^\infty \mathcal
D_{e,o}(n))q^n,$$
respectively. This proves the identity (\ref{f26}). The proof of identity (\ref{f27}) is similar.
\qed
By multiplying (\ref{f29}) with $2$ and then subtracting (\ref{f28}) we get that
\begin{corollary}\label{24} \begin{align*}
 &(-q;q)_\infty \big(\phi(-q)-2\psi(-q)-1\big)\\
=&
\sum_{n=1}^\infty(\mathcal D_{o,o}(n)+\mathcal D_{e,o}(n)-\mathcal
D_{e,e}(n)-\mathcal D_{o,e}(n))q^n\\
=&\sum_{n= 1}^\infty (-1)^{n-1} q^n(-q^{n+1})_\infty.
\end{align*}
\end{corollary}
\noindent\textbf{Proof.} To get the last identity, observe that in both sides, the coefficient of $q^n$ is the number of partitions of $n$ with the
smallest part odd minus the number of partitions of $n$ with the smallest part even.\qed
\section{Several Applications of The Involution}

The involution we use in this section is essentially the same with
previous $\alpha$, but strictly, we will denote is as $\alpha'$ to
indicate that there are minor differences.

We first prove two theorems in \mdseries{\cite{Andrews6}}
and \mdseries{\cite{GeorgeAndrews2}}. They serve as
responses to Andrews' questions on finding combinatorial
interpretations of these problems.

Recall the definition of Frobenius symbol (see \cite{Andrews1984}). The Frobenius symbol is two rows of decreasing,
non-negative integers of equal length, which is often used as one
way of representation of a partition. For example, a partition $\lambda=(7, 7,
6, 4, 4, 2, 2)$ is given, with the following Young
diagram representation. (see Figure 3)
\begin{figure}[h,t]
\setlength{\unitlength}{0.5mm}
\begin{center}
\begin{picture}(110,90)

\thicklines
            \put(0,70){\line(1,0){10}}
            \put(0,70){\line(0,-1){10}}
            \put(10,60){\line(1,0){10}}
            \put(10,60){\line(-1,0){10}}
            \put(10,60){\line(0,1){10}}
            \put(10,60){\line(0,-1){10}}
            \put(20,50){\line(1,0){10}}
            \put(20,50){\line(-1,0){10}}
            \put(20,50){\line(0,1){10}}
            \put(20,50){\line(0,-1){10}}
            \put(30,40){\line(1,0){10}}
            \put(30,40){\line(-1,0){10}}
            \put(30,40){\line(0,1){10}}
            \put(30,40){\line(0,-1){10}}
            \put(40,30){\line(0,1){10}}
            \put(40,30){\line(-1,0){10}}
                        \thinlines

    \put(0,70){\line(1,0){70}}
    \put(0,60){\line(1,0){70}}
    \put(0,50){\line(1,0){70}}
    \put(0,40){\line(1,0){60}}
    \put(0,30){\line(1,0){40}}
    \put(0,20){\line(1,0){40}}
    \put(0,10){\line(1,0){20}}
    \put(0,0){\line(1,0){20}}

        \put(0,0){\line(0,1){70}}
        \put(10,0){\line(0,1){70}}
        \put(20,0){\line(0,1){70}}
        \put(30,20){\line(0,1){50}}
        \put(40,20){\line(0,1){50}}
        \put(50,40){\line(0,1){30}}
        \put(60,40){\line(0,1){30}}
        \put(70,50){\line(0,1){20}}
        \end{picture}
\caption{ Young diagram of the partition $(7,7,6,4,4,2,2)$}
\end{center}
\end{figure}
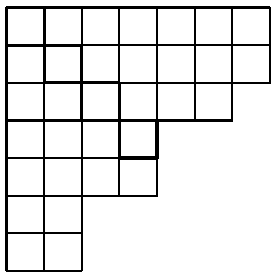

One counts the number of squares in the $4$
rows to the right and up of the diagonal, and the number of squares in the 4 columns to the
left and below the diagonal, getting the Frobenius symbol as

$${{6\ \ 5\ \ 3\ \ 0}\choose {6\ \ 5\ \ 2\ \ 1}}$$

\begin{theorem}\emph{([7, Theorem 4])}\label{31} The number of partitons of $n$
whose Frobenius symbol has no $0$ on the top row equals the number
of partitions of $n$ in which the smallest number that is not a
summand is odd.
\end{theorem}
\vskip 0.6cm
In order to prove Theorem \ref{31}, we first reduce it to proving
the following identity:
\begin{lemma}\emph{([7, Lemma 12])}\label{32}
$$\sum_{n=0}^\infty\frac{q^{n^2+n}}{(q)_n^2}=\frac{1}{(q)_\infty}\sum_{n=0}^\infty(-1)^nq^{n(n+1)/2}.$$
\end{lemma}
\vskip 0.6cm

The left hand side is the generating function for partitions of $n$
whose Frobenius symbol has no $0$ in the top row. To see this, we
draw a $n\times (n+1)$ rectangle, attach a partition with the
largest part at most $n$ under the rectangle and then attach a
partition with the largest part at most $n$ to the right of the
rectangle.

At the same time, the right hand side is actually generating
function for the number of partitions of $n$ in which the smallest
number that is not a summand is odd, since:

$$\frac{1}{(q)_\infty}\sum_{n=0}^\infty(-1)^nq^{n(n+1)/2}=
\frac{1}{(q)_\infty}\sum_{n=0}^\infty q^{n(2n+1)}(1-q^{2n+1})
=\sum_{n=0}^\infty\frac{q^{1+2+\cdots+2n}}{\prod_{j=1 \atop j\neq
2n+1}^\infty(1-q^j)}.$$

\noindent\textbf{Proof of Theorem 3.2.} Rewrite lemma \ref{32} by
multiplying $(q)_\infty$ in both sides, we get

$$\sum_{n=0}^\infty\frac{q^{n^2+n}(q^{n+1})_\infty}{(q)_n}=\sum_{n=0}^\infty(-1)^nq^{n(n+1)/2}.$$

Like in the proof of Theorem \ref{21}, we split and rearrange every term in the left-hand side. This time
the $n^2+n$ squares can be split into two copies of Sylvester's triangle $\rho_n$. The first $\rho_n$ is glued with the partition
whose largest part is at most $n$, constructing the partition $\lambda$. The second $\rho_n$ is glued with partition whose smallest part
is at leat $n+1$, and get a larger Sylverster's triangle $\rho_{n+\ell}$ and the partition $\mu$.
Then we have a triple $(\lambda,\mu,\rho_{n+\ell})$, assigned $(-1)^\ell$, such that

1. $\lambda\in \mathcal {D}$ is a partition with $n$
distinct parts,

2. $\mu\in \mathcal {P}'$ is a partition with $\ell$ parts, while the empty parts are allowed, and

3. $\rho_{n+l}$ is a Sylvester's triangle.

Similar to the previous $\alpha$, we define an involution $\alpha'$
as follows. Compare $\lambda_1$ the first part of $\lambda$ and
$\mu_1$ the first part of $\mu$. If $\lambda_1\geq \mu_1$, we remove
the first part of $\lambda$ and attach it to $\mu$. This move
changes the sign since it add $1$ to $\ell$. If $\lambda_1 < \mu_1$,
we remove the first part of $\mu$ and attach it to $\lambda$. This
move also changes the sign since it subtract $1$ from $\ell$.

 All are canceled except those that with $\lambda$ empty
partition $\emptyset$, where we have $n=0$. What have been left are
$(-1)^\ell \rho_\ell$, thus the conclusion.\qed \vskip 0.6cm

\noindent\textbf{Definition.} ([6]) A
partition with initial $k$-repetitions is a partition in which if
any $j$ appears at least $k$ times as a part then each positive
integer less than $j$ appears at least $k$ times as a part.

In a partition, one part is called a distinct part if it only
appears once. Let $D_e(m, n)$ (resp. $D_o(m, n)$) denote the number
of partitions of $n$ with initial $2$-repetitions, with $m$
different parts and an even (resp. odd) number of distinct parts.

\begin{theorem}\emph{[6, Theorem 2]}\label{33}

$$D_e(m,n)-D_o(m,n)=\begin{cases} (-1)^j, \text{ if } m=j,n=j(j+1)/2; \\
0, \text{ otherwise}.
\end{cases}$$
\end{theorem}
\vskip 0.6cm

\noindent\textbf{Proof.} We first write down the bivariate
generating function as follows:
\begin{align*}
    &\sum_{n,m\geq 0} (D_e(m,n)-D_o(m,n))z^mq^n \\
   =&\sum_{n=0}^{\infty}\frac{z^nq^{1\cdot 2+2\cdot 2+\cdots +n\cdot
2}}{(q)_n}\prod_{j=n+1}^{\infty}(1-zq^j)\\
   =& \sum_{n=0}^\infty \frac{z^nq^{n^2+n}(zq^{n+1})_\infty}{(q)_n}.
\end{align*}

Then we reduce the proof to the following identity:

\begin{equation}\label{f31}
\sum_{n=0}^\infty
\frac{z^nq^{n^2+n}(zq^{n+1})_\infty}{(q)_n}=\sum_{j=0}^\infty
(-1)^jz^jq^{\frac{j(j+1)}{2}}.
\end{equation}

Now we apply the involution $\alpha'$, while inserting a new variable $z$. In the context of the proof of Theorem \ref{32},
we use $z^{n+\ell}$ to denote the Sylvester's triangle $\rho_{n+\ell}$. Since the application does not change the value of
$n+\ell$, the proof is finished. \qed

In (\ref{f31}), let $q\to q^2$ and $z\to q^{-1}$, and we have

$$(q;q^2)_\infty\sum_{n=0}^\infty\frac{q^{2n^2+n}}{(q)_{2n}}=
\sum_{j=0}^\infty (-1)^jq^{j^2}.$$

or, by the classical identity of Euler, $(q;q^2)=1/(-q;q)_\infty$,
we have:

$$\sum_{n=0}^\infty\frac{q^{1+2+\cdots +2n}}{(q)_{2n}}=(-q)_\infty
\sum_{j=0}^\infty (-1)^jq^{j^2}.$$

Let $Q_E(n)$ denote the number of partitions of $n$ into even number
of distinct parts; $Q(n)$ denote the number of partitions of $n$ into
distinct parts. So we have proved the following identity
combinatorially.

\begin{theorem}\label{34}
\begin{equation*}
Q_E(n)=Q(n)-Q(n-1^2)+Q(n-2^2)-Q(n-3^2)+\cdots.
\end{equation*}
\end{theorem}
\noindent \noindent\textbf{Remark.} Interestingly enough, for
unrestricted partitions, the same conclusion holds. Compare the above
proof with the combinatorial proof of the following identity in \cite{Yee}:

$$p_E(n)=p(n)-p(n-1^2)+p(n-2^2)-p(n-3^2)+\cdots,$$

where $p_E(n)$ and $p(n)$ denote the number of partitions into an
even number of parts and the number of partitions, respectively.

Using the same involution of Theorem \ref{31}, we can also prove the
following two identities. We omit the proofs.

\begin{corollary}\label{35}
\begin{align*}\sum_{n=1}^\infty\frac{q^{n^2}}{(q)_{n-1}(q)_n}=&\frac{1}{(q)_\infty}
\sum_{n=1}^\infty(-1)^{n-1}q^{n(n+1)/2},
\end{align*}
or, more generally,
\begin{align*}\sum_{n=1}^\infty\frac{z^nq^{n^2}(zq^{n+1})_\infty}{(q)_{n-1}}=
\sum_{n=1}^\infty(-1)^{n-1}z^nq^{n(n+1)/2}.
\end{align*}

\end{corollary}

\begin{corollary}\label{36}
\begin{equation}\sum_{n=0}^\infty\frac{z^nq^{n^2}}{(zq)_n(q)_n}=\frac{1}{(zq)_\infty}.
\end{equation}
\end{corollary}

\noindent \textbf{Remark.} Corollary \ref{36} can be proved using
Durfee square (see \cite{GeorgeAndrews1}). Here we have a new
combinatorial proof, though a little more complicated than the
standard one.

The next application of the involution is an identity involving
another Ramanujan's third order mock theta function:
$f(q)=\sum_{n=0}^\infty \frac{q^{n^2}}{(-q)_n^2}$. Fine \cite{Fine2}
derived the following identity applying some transformation
formulas.

\begin{theorem}[\mdseries{[15, pp. 56]}]\label{37}
\begin{equation*}f(q)=1+\frac{1}{(-q;q)_\infty}\sum_{n\geq
1}(-1)^{n-1}q^n(-q^{n+1})_\infty.
\end{equation*}
\end{theorem}

\noindent \textbf{Remark.} Combinatorially, $f(q)-1 = \sum_{
n=1}^\infty (N_e(n) - N_o(n))q^n$, where $N_e(n)$ (respectively,
$N_o(n)$) is the number of partitions of $n$ with even
(respectively, odd) rank.

On the other hand, $\sum_{k\geq
1}(-1)^{k-1}q^k(-q^{k+1})_\infty=\sum_{ n=1}^\infty (L_o(n) -
L_e(n))q^n$, where $L_e(n)$ (respectively, $L_o(n)$) is the number of
partitions of $n$ into distinct parts with the smallest part even
(respectively, odd). So the above identity relates these two
enumeration problems.

Differently with previously involutions, we will compare the
smallest parts of $\mu$ and $\lambda$ along the way. All other
procedures are similar.

\noindent\textbf{Proof.} As before, we rewrite the identity as

\begin{equation}\label{f33}
\sum_{n=1}^\infty\frac{q^{n^2}(-q^{n+1})_\infty}{(-q)_n}=\sum_{n\geq
1}(-1)^{n-1}q^n(-q^{n+1})_\infty.
\end{equation}

Again, we split the $n^2$ squares into a Sylvester's triangle $\rho_n$ and another
Sylverster's triangle $\rho_{n_1}$. We glue $\rho_{n-1}$ with partition whose largest part
is at most $n$, and transpose to its conjugate to get a partition $\lambda$. Then we glue $\rho_n$ with the partition
into distinct parts whose smallest part is at least $n+1$, split it to get a larger Sylvester's triangle
$\rho_{n+\ell}$ and a partition $\mu$. Thus, we have a triple
$(\lambda,\mu,\rho_{n+\ell})$ assigned $(-1)^{k-n+1}$, such that

1. $\lambda\in \mathcal {D}'$ is a partition with $n$
distinct parts (empty part may be included). Set the largest part of
$\lambda$ as $k$.

2. $\mu\in \mathcal {P}'$ is a partition with $\ell$ parts (empty
parts may be included).

3. $\rho_{n+\ell}$ is a Sylvester's triangle.

Similar to the previous involution, we define an involution
$\alpha''$ as follows. Compare $\lambda_{small}$ the smallest part of
$\lambda$ and $\mu_{small}$ the smallest part of $\mu$ (both parts could be
empty parts). If $\lambda_{small}\leq \mu_{small}$, we remove the
smallest part of $\lambda$ and attach it to $\mu$. This move changes
the sign since it substracts $1$ from $n$. If $\lambda_{small} >
\mu_{small}$, we remove the least part of $\mu$ and attach it to
$\lambda$. This move also changes the sign since it add $1$ from
$n$.

All are canceled except those $(\lambda,\mu)$ with sign $(-1)^{k-1}$
satisfying that $\lambda$ only cantains one part $k-1$ and $\mu$
contains $\ell$ nonempty parts, each $\geq k-1$. Then attaching
this $k-1$ under $\mu$, aside with the Sylvester's triangle
$\rho_{l+1}$, we get the right-hand side of the identity.\qed \vskip
0.6cm

Combining both the above theorem and Corollary \ref{24}, we get the
following well known relation, which is first derived by Ramanujan:
\begin{corollary}\label{38}
$$ \phi(-q)-2\psi(-q) = f(q).$$
\end{corollary}

Just like Corollary \ref{35} and Corollary \ref{36}, we get two
identities of Theorem \ref{37}. The proof is also similar to that of
Theorem \ref{37}, except that this time, the partition $\lambda$ is not allowed
to have empty parts, which left more terms in the right-hand side.

\begin{corollary}\label{39}
\begin{equation}\label{f34}\sum_{n=1}^\infty\frac{q^{n^2+n}}{(-q)_n^2}=\frac{1}{(-q)_\infty} \sum_{k\geq
1} \sum_{n\geq 1} (-1)^{n-1} q^{\frac{k(k-1)}{2}} q^{k+n}
(-q^{k+n+1})_\infty.
\end{equation}
\begin{equation}\label{f35}
\sum_{n=1}^\infty\frac{q^{n^2}}{(-q)_{n-1}(-q)_n}=\frac{1}{(-q)_\infty}\sum_{k\geq
1}\sum_{n\geq
1}(-1)^{n-1}q^{\frac{k(k+1)}{2}}q^{k+n+1}(-q^{k+n+2})_\infty.
\end{equation}
\end{corollary}

As in Theorem \ref{33}, we can generalize (\ref{f33}), (\ref{f34})
and (\ref{f35}) by inserting a new variable $z$ into the identities
to track different values, getting the following:

\noindent \textbf{ Corollary 3.9*}
\begin{equation}
\sum_{n=1}^\infty\frac{z^nq^{n^2}(-zq^{n+1})_\infty}{(-q)_n}=\sum_{k\geq
1}(-1)^{k-1}zq^k(-zq^{k+1})_\infty.
\end{equation}

\begin{equation}\label{f37}
\sum_{n=1}^\infty\frac{z^nq^{n^2+n}(-zq^{n+1})_\infty}{(-q)_n}=\sum_{k\geq
1}\sum_{n\geq
1}(-1)^{n-1}z^kq^{\frac{k(k-1)}{2}}q^{k+n}(-zq^{n+k+1})_\infty.
\end{equation}

\begin{equation}
\sum_{n=1}^\infty\frac{z^nq^{n^2}(-zq^{n+1})_\infty}{(-q)_{n-1}}=\sum_{k\geq
1}\sum_{n\geq 1} (-1)^{n-1} z^{k+1} q^{\frac{k(k+1)}{2}}
q^{k+n+1}(-zq^{n+k+2})_\infty.
\end{equation}

\noindent\textbf{Proof.} For all three identities, we use $z^{n+\ell}$ to denote the Sylvester's triangle
$\rho_{n+\ell}$, which, under the application of the involution $\alpha''$, does not change. All the other procedures
follow the proof of Theorem \ref{37}.
\qed

The above three identities can be translated into the language of
partitions. We only do this for identity (3.7), stating it as a theorem,
which, as will be seen, is in the spirit of Andrews' ``partition
with initial repetitions". For more details, see \cite{Andrews6}.

Let $I_e(m, n)$ (resp. $I_o(m, n)$) denote the number of partitions
of $n$ with initial $2$-repetitions, with $m$ different parts and an
even (resp. odd) number of repeated parts. We say a partition into
distinct parts is with a initial Sylvester's triangle $\rho_k, k\geq 1$ when
it contains $1,2,\cdots,k$ and does not contain $k+1$ and it is not
$\rho_k$ itself. Let $S_e(m,n)$ (resp. $S_o(m, n)$) denote the
number of partitions of $n$ with initial sylvester's triangle, with
$m$ different parts, and the first gap (i.e. the difference between
neighboring parts which is larger than one) of parts is even (resp.
odd).

\begin{theorem}

$$I_e(m, n)-I_o(m, n)=S_o(m, n)-S_e(m,n).$$

\end{theorem}
\vskip0.6cm

\noindent\textbf{Proof.} Since in the right-hand side of (\ref{f37})
the coefficient of $z^mq^n$ is $S_o(m, n)-S_e(m,n)$. To finish the
proof, we only need to rewrite the left-hand side of (\ref{f37}) as
follows:

\begin{align*}
  &\sum_{n=1}^\infty\frac{z^nq^{n^2+n}(-zq^{n+1})_\infty}{(-q)_n}\\
 =&\sum_{n=1}^{\infty}\frac{z^nq^{1\cdot 2+2\cdot 2+\cdots +n\cdot
2}}{(-q)_n}\prod_{j=n+1}^{\infty}(1+zq^j)\\
 =& \sum_{n,m\geq 1} (I_e(m,n)-I_o(m,n))z^mq^n.
 \end{align*}
 \qed

\section{Flushed Partitions, Concave Compositions and Proper Partitions}

We first prove two identities:
\begin{lemma}
\begin{align}\label{41}
  \sum_{n=1}^\infty q^{n(3n-1)/2}(1-q^n)=\sum_{n=1}^\infty\frac{(-1)^{n-1}
q^{\frac{n(n+1)}{2}}}{(1+q)\cdots(1+q^n)}.\\
\sum_{n=1}^\infty q^{n(3n-1)/2}(1-q^n)=\sum_{n=1}^\infty
\frac{q^{n(2n-1)}}{(-q)_{2n}}.
\end{align}
\end{lemma}

\noindent\textbf{Proof.} For the second identity, let  $P_o(\mathcal
{D} _n)$ (resp. $P_e(\mathcal {D} _n)$) denote the number of
partitions $\lambda$ of $n$ with distinct parts and the largest part
is odd (resp. even). The following identity is due to Fine (see
\cite{Fine}) which can be reached by Franklin's well known
involution.
\begin{equation}\sum_{n=1}^\infty q^{n(3n-1)/2}(1-q^n)=
\sum_{n=1}^\infty(P_o(\mathcal {D}_n)-P_e(\mathcal {D}_n))q^n.
\end{equation}

To prove the first identity, it is sufficient to prove the following
\begin{equation}\label{9}
\sum_{n=1}^\infty(P_o(\mathcal {D}_n)-P_e(\mathcal
{D}_n))q^n=\sum_{n=1}^\infty\frac{(-1)^{n-1}q^{\frac{n(n+1)}{2}}}{(1+q)\cdots(1+q^n)}.
\end{equation}

Identity (\ref{9}) can be shown as follows. Given a partition with
distinct $n$ parts $\lambda\in \mathcal D$, we split the partition
into a Sylvester's triangle $\rho_n=(n,n-1,\cdots,1)$ and a
 partition $\mu\in \mathcal P$. The conjugate of $\mu$ is
 a partition, $\mu^*$, with the largest part at most $n$ and with the number
 of parts $r$. Note that $\lambda_1$, the largest part of $\lambda$ is $n+r$.
 So in the left-hand side of (\ref{9}), the above $\lambda$ is
 actually assigned $(-1)^{n-1+r}$.

 The right hand side of (\ref{9}), at the same time, is the generating function
 of a pair of partitions $(\rho_n,\mu^*)$, where $\rho_n$ is
 Sylvester's triangle and $\mu_*$ has largest part no more than $n$
 and with r parts. Each term is assigned $(-1)^{n-1+r}$, thus
 the conclusion.

 The proof of the second identity is in the same fashion of (4.4),
  while the parity of the Sylvester's triangle
 should be noted. Again, for a partition of $n$ with
 distinct parts $\lambda\in \mathcal D_n$ assigned $(-1)^{\lambda_1}$, we split it into a
 Sylvester's triangle $\rho_{2n-1}=(2n-1,2n-2,\cdots,1)$ and a
 partition $\mu$ with at most $2n$ parts, assigned
 $(-1)^{\mu_1}$, where $\mu_1$ is the largest part of $\mu$.
\qed \vskip 0.6cm

\noindent\textbf{Definition:} Durfee Symbol was introduced by
Andrews in \cite{Andrews7}, which is defined as follows: using
Corollory 3.6, we can represent an unrestricted partition as its
Durfee square of size $n$ and two partitions with the largest part
at most $n$. We then denote the partition with two rows of integers,
the top row listing the parts of the conjugate of the partition to
the right of the Durfee square and the bottom row lising the parts
of the partition under the Durfee square. We also write down a
subscript $n$ in the end to denote the size of the Durfee square.
Take the partition $\lambda=(11,11,11,9,7,5,5,4,4,3)$ for
example, whose Young diagram is depicted in Figure 4, The Durfee
symbol representation of $\lambda$ is as follows:

$${{5\ \ 5\ \ 4\ \ 4\ \ 3\ \ 3}\choose {5\ \ 5\ \ 4\ \ 4\ \ 3\ \ 0}}_5$$

Suppose the Durfee square of a partition is of size $n$, we call a
partition \emph{proper} if its Durfee symbol has the same
number of $n$'s in both the top and bottom rows. All other
partitions are \emph{improper partitions}. The number of proper
partitions of $n$ is denoted as $PR(n)$. The number of improper
partitions of $n$ is denoted as $IMPR(n)$. A typical example of
a proper partition, $\lambda=(11,11,11,9,7,5,5,4,4,3)$, has its Young
diagram in Figure 4. In its Durfee symbol, both the top and bottom rows have two $5$'s.

\begin{figure}[h,t]
\setlength{\unitlength}{0.4mm}
\begin{center}
\begin{picture}(110,100) \thicklines
\put(0,100){\line(1,0){50}} \put(0,50){\line(1,0){50}}
\put(0,50){\line(0,1){50}} \put(50,50){\line(0,1){50}} \thinlines

\put(0,100){\line(1,0){110}} \put(0,90){\line(1,0){110}}
\put(0,80){\line(1,0){110}} \put(0,70){\line(1,0){110}}
\put(0,60){\line(1,0){90}} \put(0,50){\line(1,0){70}}
\put(0,40){\line(1,0){50}} \put(0,30){\line(1,0){50}}
\put(0,20){\line(1,0){40}} \put(0,10){\line(1,0){40}}
\put(0,0){\line(1,0){30}}

\put(0,100){\line(0,-1){100}}\put(10,100){\line(0,-1){100}}
\put(10,30){\line(0,-1){30}}
\put(20,100){\line(0,-1){100}}\put(20,30){\line(0,-1){30}}
\put(30,100){\line(0,-1){100}}\put(30,30){\line(0,-1){30}}
\put(40,100){\line(0,-1){90}}
\put(40,30){\line(0,-1){20}}\put(50,100){\line(0,-1){70}}
\put(60,100){\line(0,-1){50}} \put(70,100){\line(0,-1){50}}
\put(80,100){\line(0,-1){40}} \put(90,100){\line(0,-1){40}}
\put(100,100){\line(0,-1){30}} \put(110,100){\line(0,-1){30}}

\end{picture}
\caption{$\lambda=(11,11,11,9,7,5,5,4,4,3)$}

\end{center}
\end{figure}
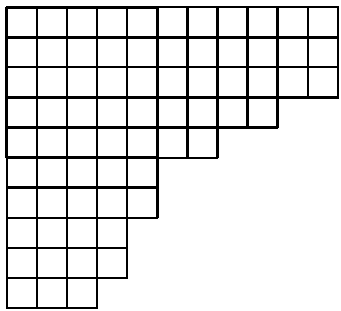

Obviously, the generating function of proper partitions is
$$\sum_{n=0}^\infty PR(n)q^n=1+\sum_{n=1}^\infty\frac{q^{n^2}}{(q)_{n-1}^2(1-q^{2n})}.$$

As we stated in the introduction, a partition $\lambda$ is called {\it
flushed} if the smallest part to appear an even number of times is
even. The unflushed partition are, of course, those partitions in
which the smallest part to appear an even number of times is odd. We
denote the number of flushed partitions of $n$ as $F(n)$.

In a flushed partition, suppose $1,2,\cdots,2i-1$ all appear an odd
number of times, and $2i$ appear an even number of times (zero times
included). We extract $1+2+\cdots+2i-1$ and left $1,2\cdots, 2i$ all
appearing an even number of times. So one easily writes down the
generating function of $F(n)$ as follows:

\begin{equation}\label{f45}
\sum_{n\geq 1} F(n) q^n= \sum_{n=1}^\infty
\frac{q^{n(2n-1)}}{(q^2;q^2)_{2n}}\frac{1}{(q^{2n+1})_\infty}=\sum_{n=1}^\infty\frac{q^{n(3n-1)/2}(1-q^n)}{(q)_\infty},
\end{equation}
where the second equality is by (4.2).

Relating Lemma \ref{41} and Theorem \ref{21} we get the following:

\begin{theorem} [Generating Functions for
Flushed Partitions]\label{42} The number of flushed partitions of
$n$ is equal to the number of proper partitions of $n$. Thus, they
share the same generating function,
$$\sum_{n\geq 1} F(n) q^n=\sum_{n=1}^\infty PR(n)q^n=
\frac{1}{(q)_\infty}\sum_{n=1}^\infty q^{n(3n-1)/2}(1-q^n).$$
\end{theorem}
\vskip 0.6cm

\noindent \textbf{Remark}: In \cite{Dyson} Dyson defined the rank of
a partition as the largest part minus the number of parts. We denote
the number of partitions of $n$ with rank m by $N(m, n)$. Then in
\cite{Atkin} Atkin and Swinnerton-Dyer derived the generating
function of $N(0,n)$ as follows:

$$\sum_{n=1}^\infty N(0,n)q^n=\frac{1}{(q)_\infty}\sum_{n=1}^\infty
(-1)^{n-1}q^{n(3n-1)/2}(1-q^n).$$

It would be interesting to compare this with Theorem \ref{42}.
Notice that there is some kind of symmetry between proper partitions
and unrestricted partitions with rank $0$. Since by the classical
pentagonal number theorem, we have $(q)_\infty=
\sum_{n=0}^\infty (-1)^{n}q^{n(3n-1)/2}(1+q^n)$, so, the above identity and
theorem \ref{42} reveal the relations of pentagonal number theorem
with three different variations of signs.

\begin{theorem} [Generating Functions for
Concave Compositions of Even Length]\label{43} The number of concave
compositions of even length of $n$ is equal to the number of
improper partitions of $n$ and thus is equal to the number of
unflushed partitions of $n$. So, they share the same generating
function,
$$
\sum_{n=0}^\infty ce(n)q^n=\sum_{n=0}^\infty
IMPR(n)q^n=\frac{1}{(q)_\infty}\Big(1-\sum_{n=1}^\infty
q^{n(3n-1)/2}(1-q^n)\Big).$$
\end{theorem}
\vskip 0.6cm

\noindent \textbf{Remark}: Theorem \ref{11} is the combination of
theorem \ref{42} and theorem \ref{43}. \\

\noindent{\textbf{Proof}}. We will construct a bijection to finish
the proof. Define a map $\phi$ from concave compositions of even
length of $n$ to unrestricted partitions of $n$ as follows. First take
a concave composition $\mathcal{C}:$
 $$a_1>a_2>\cdots>a_m=b_m<b_{m-1}<\cdots<b_1.$$
The partition  $\phi(\mathcal{C})$ will depend on four cases.

1. When $a_m=b_m=0, a_{m-1}=b_{m-1}$, let the first row of the
Young diagram have $a_1$ squares. Then draw $b_1$ squares under
the square $(1,1)$; Begin with the position $(2,2)$, $a_2$ squares
are put in the second row. Then we put $b_2$ squares under the
position $(2,2)$. Continuing in this fashion, we get a diagram. This
diagram represent a partition of $n$, which is denoted as $\phi(\mathcal C)$.
Now we claim that the resulting partition is a improper partition.
In fact, since $a_m=b_m=0$, one knows that the Durfee square of
$\phi(\mathcal C)$ is with size $m-1$. Thus, in its Durfee symbol, in the
top row, the number of
$m-1$'s is $a_{m-1}-1$; in the bottom row, the number of $m-1$'s is $b_{m-1}$.
For example, if the concave composition is $\mathcal{C}:2>1>0=0<1<2$, we
have $\phi(\mathcal{C})=(2,2,2)$, as illustrated in Figure 5.

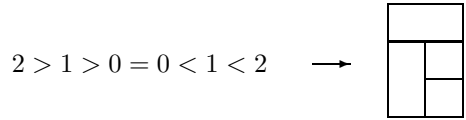
\begin{figure}[h,t]
\setlength{\unitlength}{0.5mm}
\begin{center}
\begin{picture}(120,50)

\put(0,12){$2>1>0=0<1<2$}\put(80,14){\vector(1,0){10}}

\put(100,0){\line(1,0){20}} \put(110,10){\line(1,0){10}}
\put(100,20){\line(1,0){20}} \put(100,30){\line(1,0){20}}

\put(100,0){\line(0,1){30}}\put(110,0){\line(0,1){20}}\put(120,0){\line(0,1){30}}

\end{picture}\caption{concave compositions with two second smallest equal parts}
\end{center}
\end{figure}

2. When $a_m=b_m\neq0$, we do symmetrically as in case one. let the first column of the Young diagram have $b_1$ squares.
Then draw $a_1$ squares to the right of the square $(1,1)$; Begin
with the position $(2,2)$, $b_2$ squares are put in the second
column. Then we put $a_2$ squares to the right of the position
$(2,2)$. Continuing in this fashion, we get a Young diagram. This diagram
represent a partition, which, denote by $\phi(\mathcal C)$, is improper.
To see this, we first observe that the Durfee square of $\phi(\mathcal C)$
is of size $m$. In the top row of the Durfee symbol, the number of $m$'s is
$a_m$, while in the bottom row, the number of $m$'s is $b_m-1$.
For
example, if the concave composition is $\mathcal{C}:2>1=1<2$, we
have $\phi(\mathcal{C})=(3,3)$ as illustrated in Figure 6.

\begin{figure}[h,t]
\setlength{\unitlength}{0.5mm}
\begin{center}
\begin{picture}(120,50)

\put(0,12){$2>1=1<2$}\put(60,14){\vector(1,0){10}}

\put(80,0){\line(1,0){30}} \put(80,20){\line(1,0){30}}
\put(80,0){\line(0,1){20}} \put(90,0){\line(0,1){20}}

\put(90,10){\line(1,0){20}}\put(100,0){\line(0,1){10}}\put(110,0){\line(0,1){20}}

\end{picture}\caption{concave compositions without zeros}
\end{center}
\end{figure}
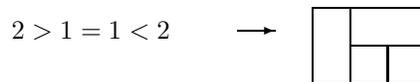

3. When $a_m=b_m=0, a_{m-1}\neq b_{m-1}$, if $a_{m-1}> b_{m-1}$, we
follow the procedure of case two. We get a diagram, representing a partition of
$n$, whose Durfee square has size $m-1$. We denote this partition as $\phi(\mathcal C)$
and observe that, in its Durfee symbol, the number of $m-1$'s in the top row is $a_{m-1}$, while
the number of $m-1$'s in the bottom row is $b_{m-1}-1$. Since $a_{m-1}-b_{m-1}+1\geq 2$, we
have that $\phi (\mathcal C)$ is improper.

For example, if the concave composition is
$\mathcal{C}:3>2>0=0<1<2$, we have $\phi(\mathcal{C})=(4,4)$, as
illustrated in Figure 7.

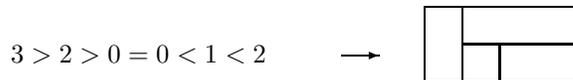
\begin{figure}[h,t]
\setlength{\unitlength}{0.5mm}
\begin{center}
\begin{picture}(220,50)

\put(20,15){$3>2>0=0<1<2$}

\put(108,17){\vector(1,0){10}}

 \put(150,10){\line(1,0){20}}
\put(140,20){\line(1,0){30}}
\put(130,30){\line(1,0){40}}

\put(130,10){\line(0,1){20}}

\put(140,10){\line(0,1){20}}
\put(150,10){\line(0,1){10}}

\put(170,10){\line(0,1){20}}
\put(130,10){\line(1,0){40}}

\end{picture}\caption{the map $\phi$: $3>2>0=0<1<2 \to (4,4)$}
\end{center}
\end{figure}

4. When $a_m=b_m=0, a_{m-1}\neq b_{m-1}$, if $a_{m-1}> b_{m-1}$. We
follow the procedure of case one. We get a diagram which represents a partition of
$n$, whose Durfee square has size $m-1$. We denote this partition as $\phi(\mathcal C)$
and observe that, in its Durfee symbol, the number of $m-1$'s in the top row is $a_{m-1}-1$, while
the number of $m-1$'s in the bottom row is $b_{m-1}$. Since $b_{m-1}-a_{m-1}+1\geq 2$, we
have that $\phi (\mathcal C)$ is improper.

Now we finished the definition of the map $\phi$, which map all
concave compositions of even length of $n$ into improper partitions
of $n$. We still need to see that $\phi$ is a bijection. Given an
improper partition, $\mathcal C'$, and suppose for accuracy its
Durfee square's size is $m$. We then write down its Durfee symbol, and observe
the difference of the number of $m$'s in its top row and bottom tow. We divide it into
four cases, namely, when the difference is $1, -1,\geq 2 \text{ and } \leq -2$.
Each of these cases corresponds to one the four cases we analyzed above, then we can map them back
to the concave compositions of even length, which shows $\phi$ is indeed a bijective map.
\qed

In order to prove Sylvester's first problem, we have the following
Corollary of the above theorems.

\begin{corollary}\label{44}
The number of unflushed partitions of $n$ with $m$ parts is equal
with the number of unrestricted partitions of $n$ with $m$ parts minus the number
of partitions of $n$ with $m$ parts in
which when the Durfee square is $k\times k$, the Durfee symbol
contains no $k$'s in the top row and contains an even number of $k$'s in the bottom row.
\end{corollary}

\noindent \textbf{Proof.} Just like what we have done in the
previous two sections, we try to insert another variable $z$ into
the identities and start all the procedures of finding generating
function of unflushed partitions all over again. Then the generating
function for unflushed partitions of $n$ with $m$ parts can be
denoted as:

\begin{align*}
&\sum_{n=0}^\infty
\frac{z^{2n}q^{n(2n+1)}}{(q^2;q^2)_{2n+1}}\frac{1}{(zq^{2n+2})_\infty}\\
=&\frac{1}{(zq)_\infty}\sum_{n=0}^\infty
\frac{z^{2n}q^{n(2n+1)}}{(-zq)_{2n+1}}\\
=&\frac{1}{(zq)_\infty}\sum_{n=0}^\infty
\frac{(-1)^nz^{n}q^{n(2n+1)/2}}{(-zq)_n}\\
=&\sum_{n=0}^\infty\frac{z^nq^{n^2}}{(zq)_n(q)_n} -\sum_{n=1}^\infty
\frac{z^{n}q^{n^2}}{(q)_{n-1}(zq)_{n-1}(1-z^2q^{2n})}.
\end{align*}
the first identity is by an analysis similar to that of (4.5), the
second identity is from the proof of Lemma \ref{41}, and the third
identity is by (2.3) and (3.2).

We interpret the first and the last sums of the above identity and
get the conclusion. \qed

We conclude this paper with a combinatorial proof of Sylvester's problem.

\begin{theorem} Unflushed partitions of $n$ with odd number of
parts are equinumerous with unflushed partitions of $n$ with even
number of parts.
\end{theorem}
\noindent\textbf{Proof.} By Corollary \ref{44}, we only need to
analyze the parity of the number of unrestricted partitions of $n$
with $m$ parts, in which, when the Durfee square is $k\times k$, either
there are an odd number of $k$'s in the bottom row, or there are an even number
of $k$'s in the bottom row and at least one $k$ in the top row.

If the total number of $k$'s in the Durfee symbol is odd, say,
$2i-1$, then we can arrange these $k$'s in top and bottom rows in
$2i$ ways, by putting $0,1,\cdots,2i-1$ copies of $k$'s in the
bottom row and the rest $k$'s in the top row. Observe that the
parity of number of parts changes accordingly, that is, among these
$2i$ partitions, $i$ ones have an even number of parts and $i$ ones
have an odd number of parts.

If, however, the total number of $k$'s in the Durfee symbol is even,
say, $2i$, remember that the top row must contain at leat one $k$.
Then we arrange these $k$'s in top and bottom
rows in $2i$ ways, by putting $0,1,\cdots,2i-1$ copies of $k$'s in
the bottom row and the rest in the top row. In this case the
parity of number of parts changes accordingly, too, that is, among
these $2i$ partitions, $i$ ones have an even number of parts and $i$
ones have an odd number of parts. \qed

\section{Possible Further Works}

1. In \cite{Sylvester}, Sylvester also asked another problem, as he put it,
``2. Required to prove that the same proposition holds when any odd number is
partitioned without repetitions in every possible way." We don't know whether
or not a similar combinatorial proof can be given for this
second problem.

\noindent 2. Regarding Section 3, it seems the method in the proof of Theorem 2.1 involving the involutions
$\alpha, \alpha' \text{ and } \alpha''$ can be applied to more problems.

\noindent\textbf{Acknowledgement.} The author thanks Mrs. Xiu You for her constant support and encouragement,
and the referees for their help and valuable suggestions.

\end{document}